\documentclass{amsart}
\date{\today}
\theoremstyle{plain}
  \newtheorem{thm}{Theorem}
\theoremstyle{remark}
  \newtheorem{rem}{Remark}
  \newtheorem*{ack}{Acknowledgment}
\newcommand{\C}{{\mathbf{C}}}
\newcommand{\cM}{{\mathcal{M}}}
\begin{document}
\title[Hyperbolic three-manifolds with trivial finite type invariants]
{Hyperbolic three-manifolds \\ with trivial finite type invariants}
\author{Hitoshi Murakami}
\address{Department of Mathematics, School of Science and Engineering,
Waseda University, Ohkubo, Shinjuku, TOKYO 169-8555, JAPAN,
and
Mittag-Leffler Institute, Aurav{\"a}gen 17, S-182 62,
Djursholm, SWEDEN}
\email{hitoshi@uguisu.co.jp}
\thanks{Partially supported by Waseda University Grant
for Special Research Projects (No. 98A-623) and
Grant-in-Aid for Scientific Research (C) (No. 09640135), the Ministry of
Education, Science, Sports and Culture.}
\begin{abstract}
We construct a hyperbolic three-manifold with trivial finite type invariants
up to an arbitrarily given degree.
\end{abstract}
\maketitle
The concept of finite type invariants for integral homology three-spheres was
introduced by T.~Ohtsuki in \cite{Ohtsuki:JKNOT96} including the Casson
invariant \cite{Akbulut/McCarthy:Casson_invariant} for the first non-trivial
example.
It attracts not only mathematicians but also physicists since it is closely
related to E.~Witten's quantum invariants for three-manifolds
\cite{Witten:COMMP89} in the following way.
\par
It is proved by T.T.Q.~Le \cite{Le:AMEMS97} that the degree $d$ term of the
LMO invariant \cite{Le/Murakami/Ohtsuki:TOPOL98} is of finite type of degree
$d$ and conversely any finite type invariant of degree $d$ comes from the
degree $d$ term of the LMO invariant.
Since it is also known that the perturbative $PSU(n)$ invariant can be obtained
from the LMO invariant \cite{Ohtsuki:perturbative_invariant,Le:1998},
every coefficient of the perturbative $PSU(n)$ invariant is of finite
type of degree $d$.
(For $n=2$ case was also proved in \cite{Kricker/Spence:JKNOT97}.)
Therefore we may say that finite type invariants approximate quantum invariants
from lower degree parts.
We refer the reader to \cite{Ohtsuki:Oiwake96,Ohtsuki:1998} for more detail.
\par
It is a natural question to ask how strong finite type invariants are.
For any given integral homology three-sphere $M$ and any positive integer $k$,
Le \cite{Le:AMEMS97}, and N.~Habegger and K.~Orr \cite{Habegger/Orr:1997} gave
infinitely many homology spheres which have
the same finite type invariants as $M$ up to degree $k$ by using the LMO
invariant.
E.~Kalfagianni \cite{Kalfagianni:MATRL97} also obtained such manifolds using $\pm1$-surgery along knots.
\par
In this paper we give at least one hyperbolic three-manifold with trivial finite
type invariants up to an arbitrarily given degree by simple construction using
a hyperbolic Brunnian link.
\begin{ack}
This work was done when the author was visiting the Mittag-Leffler Institute.
He thanks the staffs for their hospitality: `{\it Tack s{\aa} mycket}'.
He also thanks Kazuo Habiro, Thang Le and Tomotada Ohtsuki for useful
comments.
He is grateful to Nathan Habegger, who sent him his preprint
with K.~Orr \cite{Habegger/Orr:1997} where he found that the technique using a
Brunnian link has been already used in the proof of Theorem 8.1,
and Effie Kalfagianni, who informed him of her result
\cite{Kalfagianni:MATRL97}.
\par
Thanks are also due to the subscribers to the mailing list `knot' (http://w3.to/oto/) run by Makoto Ozawa, where he could enjoy useful comments.
\end{ack}
\section{Preliminaries}
In this section we describe a filtration of integral homology three-spheres
(closed, oriented three-manifolds with the same integral homology groups as
the three-sphere) and define finite type invariants following T.~Ohtsuki
\cite{Ohtsuki:JKNOT96}.
\par
Let $\cM$ be the $\C$-vector space spanned by all the integral homology
three-spheres.
A link $L$ in an integral homology three-sphere $M$ is called algebraically
split if every linking number vanishes.
A unit-framed link is a link with each framing $1$ or $-1$.
For an algebraically split, unit-framed link $L$ in an integral homology
three-sphere $M$, we put
\begin{equation*}
  [M,L]=\sum_{\ell\subset L}(-1)^{\sharp(\ell)}\chi(M;\ell)\in\cM,
\end{equation*}
where the summation runs over all sublinks of $L$ including the empty link
and $L$ itself, $\sharp(\ell)$ is the number of components in $\ell$
and $\chi(M;\ell)$ is the manifold obtained from $M$ by surgery along
$\ell$ respecting the framing.
Note that since $L$ is algebraically split and unit-framed, $\chi(M;\ell)$
is again an integral homology three-sphere.
\par
Let $\cM_{d}$ be the vector subspace of $\cM$ spanned by all $[M,L]$ with
$\sharp(L)=d$.
Now we have the following filtration of $\cM$:
\begin{equation*}
  \cM=\cM_{0}\supset\cM_{1}\supset\cM_{2}\supset\dots\supset\cM_{d}
  \supset\cM_{d+1}\supset\cdots.
\end{equation*}
An invariant $T$ of integral homology three-spheres is said to be of
{\em finite type} of degree $d$ if $T\vert_{\cM_{d+1}}=0$.
Here we extend $T$ linearly over $\cM$.
\section{Main result}
In this section we state our main result and give the proof.
Main tools are hyperbolic Brunnian links due to T.~Kanenobu and hyperbolic Dehn
surgery due to W.~Thurston.
\begin{thm}
For any positive integer $k$, there exists (at least) one hyperbolic
three-manifold $M$ such that for any finite type invariant $T_d$ of degree
$d<k$, $T_d(M)=T_d(S^3)$.
Note that $M$ is not homeomorphic to $S^3$ since it is hyperbolic.
\end{thm}
\begin{proof}
From the definition of finite type invariants it is sufficient to
construct a hyperbolic three-manifold $M$ satisfying $M-S^3\in\cM_{k}$ for a
given $k$.
\par
By \cite{Kanenobu:JMATS386}, there exists a hyperbolic Brunnian link $L$
with $k$ components, that is, the exterior of $L$ is hyperbolic and each proper
sublink of $L$ is trivial.
Let $M_{1/n_1,1/n_2,\dots,1/n_k}$ be the closed three-manifold obtained from
$S^3$ by Dehn surgery along $L$ with surgery coefficient
$(1/n_1,1/n_2,\dots,1/n_k)$ for positive integers $n_i$.
Here we use the notation described in \cite{Rolfsen:1990}.
By W.~Thurston's hyperbolic Dehn surgery theorem \cite{Thurston:BULAM382},
$M_{1/n_1,1/n_2,\dots,1/n_k}$ is hyperbolic if $n_i$ are sufficiently large for
all $i$.
Let $M=M_{1/n_1,1/n_2,\dots,1/n_k}$ be a hyperbolic three-manifold obtained
in this way.
\par
We show that $M-S^3\in\cM_{k}$.
Since each component of $L$ is trivial, $M$ is obtained by Dehn surgery
along the link $L^{n_1,n_2,\dots,n_k}$ with every coefficient one, where
$L^{n_1,n_2,\dots,n_k}$ is the $(n_1+n_2+\dots+n_k)$-component link obtained
by replacing $i$th component of $L$ with $n_i$ parallels.
\par
Now one can easily see that for $d<k$ each $d$-component sublink of
$L^{n_1,n_2,\dots,n_k}$ is trivial.
Let $L_k$ be a $k$-component sublink of $L^{n_1,n_2,\dots,n_k}$.
Then since any proper sublink in $L_k$ is trivial, we have
\begin{equation}\label{eq:M_k}
  \cM_k\ni[S^3,L_k]=(-1)^{k}\left(\chi(L_k)-S^3\right),
\end{equation}
where $\chi(L_k)$ is the three-manifold obtained by Dehn surgery along
$L_k$ with every coefficient one.
Next let $L_{k+1}$ be a $(k+1)$-component sublink of $L^{n_1,n_2,\dots,n_k}$.
Then
\begin{align*}
  \cM_{k+1}\ni[S^3,L_{k+1}]
  &=(-1)^{k+1}\chi(L_{k+1})
  +(-1)^{k}\sum_{\ell_k\subset L_{k+1},\sharp(\ell_k)=k}\chi(\ell_k)
  +(-1)^{k+1}k\,S^3
  \\
  &=(-1)^{k+1}
    \left\{
      \left(\chi(L_{k+1})-S^3\right)
      -\sum_{\ell_k\subset L_{k+1},\sharp(\ell_k)=k}
         \left(\chi(\ell_k)-S^3\right)
    \right\}.
\end{align*}
But from \eqref{eq:M_k} we see that for any $k$-component sublink $\ell_k$
of $L_{k+1}$, $\chi(\ell_k)-S^3\in\cM_{k}$.
Noting that $\cM_{k}\supset\cM_{k+1}$, we see that
$\chi(L_{k+1})-S^3\in\cM_{k}$.
Continuing this argument we can show that
$M\left(=\chi(L^{n_1,n_2,\dots,n_k})\right)-S^3\in\cM_{k}$,
completing the proof.
\end{proof}
\begin{rem}
During the preparation of this paper, the author was informed that
K.~Habiro obtained a stronger result by using his clasper theory.
\end{rem}
\bibliography{mrabbrev,hitoshi}

\ifx\undefined\bysame
\newcommand{\bysame}{\leavevmode\hbox to3em{\hrulefill}\,}
\fi
\begin{thebibliography}{10}

\bibitem{Akbulut/McCarthy:Casson_invariant}
S.~Akbulut and J.D. McCarthy, {\em {C}asson's {I}nvariant for {O}riented
  {H}omology $3$-spheres - {A}n {E}xposition}, Mathematical Notes, vol.~36,
  Princeton University Press, Princeton, 1990.

\bibitem{Habegger/Orr:1997}
N.~Habegger and K.~Orr, {\em Milnor link invariants and quantum 3-manifold
  invariants}, preprint,1997.

\bibitem{Kalfagianni:MATRL97}
E.~Kalfagianni, {\em Homology spheres with the same finite type invariants of
  bounded orders}, Math. Res. Lett. {\bf 4} (1997), no.~2-3, 341--347.

\bibitem{Kanenobu:JMATS386}
T.~Kanenobu, {\em Hyperbolic links with {B}runnian properties}, J. Math. Soc.
  Japan {\bf 38} (1986), no.~2, 295--308.

\bibitem{Kricker/Spence:JKNOT97}
A.~Kricker and B.~Spence, {\em Ohtsuki's invariants are of finite type}, J.
  Knot Theory Ramifications {\bf 6} (1997), no.~4, 583--597.

\bibitem{Le:1998}
T.~T.~Q. Le, {\em On perturbative {PSU}(n) invariants of rational homology
  3-spheres}.

\bibitem{Le:AMEMS97}
\bysame, {\em An invariant of integral homology $3$-spheres which is universal
  for all finite type invariants}, Solitons, geometry, and topology: on the
  crossroad (Providence, RI), Amer. Math. Soc. Transl. Ser. 2, vol. 179, Amer.
  Math. Soc., Providence, RI, 1997, pp.~75--100.

\bibitem{Le/Murakami/Ohtsuki:TOPOL98}
T.~T.~Q. Le, J.~Murakami, and T.~Ohtsuki, {\em On a universal perturbative
  invariant of $3$-manifolds}, Topology {\bf 37} (1998), no.~3, 539--574.

\bibitem{Ohtsuki:perturbative_invariant}
T.~Ohtsuki, {\em The perturbative ${SO}(3)$ invariant of rational homology
  $3$-spheres recovers from the universal perturbative invariant}, to appear in
  Topology.

\bibitem{Ohtsuki:Oiwake96}
\bysame, {\em Combinatorial quantum method in $3$-dimensional topology},
  Lecture Notes of the workshop at Oiwake Seminar House, Waseda University,
  Karuizawa, September, 1996.

\bibitem{Ohtsuki:JKNOT96}
\bysame, {\em Finite type invariants of integral homology $3$-spheres}, J. Knot
  Theory Ramifications {\bf 5} (1996), no.~1, 101--115.

\bibitem{Ohtsuki:1998}
\bysame, {\em A filtration of the set of integral homology $3$-spheres},
  Proceedings of the International Congress of Mathematicians, Vol. II (Berlin,
  1998), 1998, pp.~473--482.

\bibitem{Rolfsen:1990}
D.~Rolfsen, {\em Knots and links}, Publish or Perish Inc., Houston, TX, 1990,
  Corrected reprint of the 1976 original.

\bibitem{Thurston:BULAM382}
W.~P. Thurston, {\em Three-dimensional manifolds, {K}leinian groups and
  hyperbolic geometry}, Bull. Amer. Math. Soc. (N.S.) {\bf 6} (1982), no.~3,
  357--381.

\bibitem{Witten:COMMP89}
E.~Witten, {\em Quantum field theory and the {J}ones polynomial}, Comm. Math.
  Phys. {\bf 121} (1989), 351--399.

\end{thebibliography}
\bibliographystyle{amsplain}
\end{document}